\documentclass[12pt,a4paper]{amsart}
\usepackage{amsmath,amssymb,amscd}
\input xy
\xyoption{all}
\newcommand{\ra}{\rightarrow}

\newcommand{\PP}{\mathbb P}

\newcommand{\cE}{\mathcal{E}}
\newcommand{\cF}{\mathcal{F}}
\newcommand{\cL}{\mathcal{L}}
\newcommand{\cU}{\mathcal{U}}
\newcommand{\cO}{\mathcal{O}}
\newcommand{\cM}{\mathcal{M}}

\theoremstyle{plain}
\newtheorem{theorem}{Theorem}[section]
\newtheorem{lemma}[theorem]{Lemma}
\newtheorem{proposition}[theorem]{Proposition}

\begin{document}
\title[Poincar\'e bundles]{On Poincar\'e bundles of vector bundles on curves}

\author{H. Lange}
\author{P. E. Newstead}

\address{H. Lange\\Mathematisches Institut\\
              Universit\"at Erlangen-N\"urnberg\\
              Bismarckstra\ss e $1\frac{ 1}{2}$\\
              D-$91054$ Erlangen\\
              Germany}
              \email{lange@mi.uni-erlangen.de}
\address{P.E. Newstead\\Department of Mathematical Sciences\\
              University of Liverpool\\
              Peach Street, Liverpool L69 7ZL, UK}
\email{newstead@liv.ac.uk}
\thanks{Both authors are members of the research group VBAC (Vector Bundles on Algebraic Curves). The second author 
         acknowledges support from EPSRC Grant No. EP/C515064, and would like to thank the Mathematisches Institut der Universit\"at 
         Erlangen-N\"urnberg for its hospitality}

\begin{abstract}
Let $M$ denote the moduli space of stable vector bundles of rank $n$ and fixed determinant of degree coprime to $n$
on a non-singular projective curve $X$ of genus $g \geq 2$.
Denote by $\cU$ a universal bundle on $X \times M$. We show that, for $x,y \in X,\; x \neq y$, the restrictions 
$\cU|\{x\} \times M$ and $\cU|\{y\} \times M$ are stable and non-isomorphic when considered as bundles on $X$.
\end{abstract}
\maketitle

\section{Introduction}

\noindent
Let $X$ be a non-singular projective curve of genus $g \geq 2$ over the field of complex numbers.
We denote by $M = M(n,L)$ the moduli space of stable vector bundles of rank $n$ with determinant $L$ of degree $d$ 
on $X$, where gcd$(n,d) = 1$. We denote by $\cU$ a universal bundle on $X \times M$. For any $x \in X$ we denote by 
$\cU_x$ the bundle $\cU|\{x\} \times M$ considered as a bundle on $M$.

In a paper of M. S. Narasimhan and S. Ramanan \cite{nr} it was shown that $\cU_x$ is a simple bundle and that 
the infinitesimal deformation map
\begin{equation} \label{eq1}
T_{X,x} \ra H^1(M, \mbox{End} (\cU_x))
\end{equation}
is bijective for all $x \in X$. In \cite[Proposition 2.4]{bbn} it is shown that $\cU_x$ is semistable with respect to the unique polarization
of $M$. In fact, $\cU_x$ is stable; since we could not locate a proof of this in the literature, we include one here.

Let $\cM$ denote the moduli space of stable bundles on $M$ having the same Hilbert polynomial as $\cU_x$. Then (\ref{eq1})
implies that the natural morphism
$$
X \ra \cM
$$ is \'etale and surjective onto a component $\cM_0$ of $\cM$.

It is stated in \cite{nr} that it can be easily deduced from the results of that paper that
the map $X \ra \cM_0$ is also injective. This would imply that the curve $X$ can be identified 
with $\cM_0$. However no proof of this fact seems to be given. There is a proof in a paper 
of A. N. Tyurin \cite[Theorem 2]{tyu}, but this seems to us to be incomplete. We offer here a proof which is in the spirit of 
\cite{tyu}. To be more precise, our main result is the following theorem.\\

\noindent
{\bf Theorem}
{\it Let $X$ be a non-singular projective curve of genus $g \geq 2$. If $x,y \in X, \; x \neq y$, then $\cU_x \not\simeq \cU_y$.}\\

Note that if $X$ is a general curve of genus $g \geq 3$ or any curve of genus 2, then
$X$ does not admit \'etale coverings $X \ra \cM_0$ of degree $>1$. 
So for such curves the theorem is immediate. For the proof 
we can therefore assume that $g \geq 3$. In fact, our proof fails for $g=2$.

In Section 2 we prove the stability of $\cU_x$. In Sections 3 and 4 we make some cohomological computations, from which a family 
of stable bundles on $X$ can be constructed. This construction is carried out in Section 5, 
where we also use the morphism to $M$ given by this family in order to prove the theorem.

\section{Stability of $\cU_x$}

\noindent
Let $X$ be a non-singular projective curve of genus $g \geq 2$. 
Let $n \geq 2$ and $d$ be integers with gcd$(n,d) = 1$.
There are uniquely determined integers $l$ and $e$ 
with $0<l<n$ and $0 \leq e<d$ such that 
\begin{equation} \label{eq2}
ld-en = 1.
\end{equation}
The bundles $\cU_x$ were shown to be semistable in \cite[Proposition 2.4]{bbn}, but the proof does not seem to imply stability 
directly, even though we know also by \cite{nr} that $\cU_x$ is simple. 

\begin{proposition} \label{propos2.1}
For all $x \in X$, the vector bundle $\cU_x$ is stable with respect to the unique polarization of $M$.
\end{proposition}

\begin{proof}
By \cite[Proposition 2.4]{bbn} the bundle $\cU_x$ is semistable. By \cite[Remark 2.9]{ram} and possibly after tensoring $\cU$
by a line bundle on $M$,
$$
c_1(\cU_x) = l \alpha,
$$
where $\alpha$ is the positive generator of $H^2(M)$. By (\ref{eq2}), $l$ and $n$ are coprime. It follows that $\cU_x$ is stable. 
\end{proof}

\section{Cohomological constructions}

\noindent
Let $l$ and $n$ be as in (\ref{eq2}).
Let $V$ be a semistable vector bundle of rank $l$ and degree $l(n-l)+e$ and $W$ a semistable  bundle of rank $n-l$ and 
degree $d-e-l(n-l)$ on $X$. Then 
$$
\deg (W^* \otimes V) = nl(n-l) -1.
$$ 
Let $q_i, \; i=1,2,$ denote the projections of $X \times X$ on the two factors, $\Delta$ the diagonal of $X \times X$ and write for brevity
$$
U = q_1^*(W^* \otimes V).
$$

\begin{lemma} \label{lem2.1}
For $ n \geq 2$ and $1 \leq i \leq n$, \\
{\em (a)} $h^0(U(-i\Delta)|\Delta) = (n+(2i-1)(g-1))l(n-l) -1$;\\
{\em (b)} $h^1(U(-i\Delta)|\Delta) = 0$.
\end{lemma}

\begin{proof}
Identifying $\Delta$ with $X$, we have $U(-i\Delta)|\Delta = W^* \otimes V \otimes K_X^i$. Since 
$$
\deg (W^* \otimes V \otimes K_X^i) = (n+(2g-2)i)l(n-l) - 1 > l(n-l)(2g-2)
$$
and $W^* \otimes V$ is semistable, (b) holds and Riemann-Roch gives (a).
\end{proof}

\begin{lemma} \label{lem2.2} For $n \geq 2$,
$$
h^1(U(-n\Delta)) = gh^0(W^* \otimes V) + l(n-l)(n-1)(g(n-1) + 1) -(n-1).
$$
\end{lemma}

\begin{proof} 
For $0 \leq i \leq n$, consider the exact sequence
\begin{equation} \label{eqn12}
0 \ra U(-(i+1)\Delta) \ra U(-i\Delta) \ra U(-i\Delta)|\Delta \ra 0
\end{equation}
on $X \times X$. For $i=0$, this sequence gives
$$
0 \ra H^1(U(-\Delta)) \ra H^1(U) \stackrel{\psi}{\ra} H^1(U|\Delta),
$$
since the restriction map $H^0(U) \ra H^0(U|\Delta)$ is an isomorphism. The map $\psi$ is surjective, since its restriction
to the K\"unneth component $H^1(W^* \otimes V) \otimes H^0(\cO) \subset H^1(U)$ is an isomorphism. Hence
$$
\begin{array}{ll}
h^1(U(-\Delta)) &= h^1(U) - h^1(U|\Delta)\\
&= h^1(W^* \otimes V) h^0(\cO) + h^0(W^* \otimes V)h^1(\cO) - h^1(W^* \otimes V)\\
&= g\cdot h^0(W^* \otimes V).
\end{array}
$$
For $1 \leq i \leq n-1$, the sequence (\ref{eqn12}) gives, by Lemma \ref{lem2.1} (b),
$$
0 \ra H^0(U(-i\Delta)|\Delta) \ra H^1(U(-(i+1)\Delta)) \ra H^1(U(-i\Delta)) \ra 0.
$$
This gives, by Lemma \ref{lem2.1} (a) and the above computation,
$$
\begin{array}{l}
h^1(U(-n\Delta)) = h^1(U(-\Delta)) + \sum_{i=1}^{n-1} h^0(U(-i\Delta)|\Delta)\\
\quad = g h^0(W^* \otimes V) + \sum_{i=1}^{n-1}((n+(2i-1)(g-1))l(n-l) -1)\\
\quad = g h^0(W^* \otimes V) + l(n-l)(n-1)(g(n-1) + 1) - (n-1).
\end{array}
$$
\end{proof}

\begin{lemma} \label{lem2.3}
Let $n \geq 2$ and $x \in X$. Then, except in the case when $n=2$ and $W^* \otimes V \simeq \cO(x)$,
$$
h^1(U(-n\Delta-X\times\{x\}) = h^1(U(-\Delta - X \times \{x\})) + l(n-l)(n-1)^2g -(n-1).
$$
\end{lemma}

\begin{proof}
For $1 \leq i \leq n-1$ consider the exact sequence
$$
\begin{array}{ll}
0 \ra U(-(i+1)\Delta - X \times \{x\}) \ra & U(-i\Delta-X \times \{x\})\\ 
   &\ra U(-i\Delta-X \times \{x\})|\Delta \ra 0
\end{array}
$$
on $X \times X$.
Identifying $\Delta$ with $X$, we have 
$$
U(-i\Delta-X \times \{x\})|\Delta \simeq K_X^i \otimes W^* \otimes V(-x).
$$
If either $i \geq 2$ or $ n \geq 3$,
$$
\deg(K_X^i \otimes W^* \otimes V(-x)) > l(n-l)(2g-2).
$$
So semistability implies 
\begin{equation} \label{eqn2}
h^1(K_X^i \otimes W^* \otimes V(-x)) = 0.
\end{equation}
If $n=2$ and $i=1$, then $W^* \otimes V$ has rank 1 and 
$$
\deg(K_X \otimes W^* \otimes V(-x)) = 2g-2.
$$
So (\ref{eqn2}) is still true, 
unless $W^* \otimes V \simeq \cO(x)$.

Now Riemann-Roch implies
$$
h^0(K_X^i \otimes W^* \otimes V(-x)) = ((2g-2)i + n-g)l(n-l) -1.
$$
Hence applying the above sequence $n-1$ times, we get

$$
\begin{array}{l}
h^1(U(-n\Delta-X \times \{x\}) = \\
\quad = h^1(U(-\Delta - X \times \{x\})) + \sum_{i=1}^{n-1} h^0(K_X^i \otimes W^* \otimes V(-x))\\
\quad = h^1(U(-\Delta -X \times \{x\})) + \sum_{i=1}^{n-1} \{((2g-2)i +n-g)l(n-l) -1\}\\
\quad = h^1(U(-\Delta -X \times \{x\})) + l(n-l)(n-1)^2g -(n-1).
\end{array}
$$
\end{proof}

Now suppose $(V,W)$ is a general pair of bundles on $X$ with the given ranks and degrees. Here by ``general'' we mean that
the theorem of Hirschowitz (see \cite{hir}) is true, which says that either $H^0(W^* \otimes V) = 0$ or $H^1(W^* \otimes V) = 0$.

\begin{proposition} \label{prop2.4}
For $n \geq 3$, $g \geq 3$ and $(V,W)$ general, there is a 2-dimensional vector subspace $T_0 \subset H^1(U(-n\Delta))$ 
such that the restriction map
\begin{equation} \label{eqn3}
H^1(U(-n\Delta)) \ra H^1(W^* \otimes V(-nx))
\end{equation}
is injective on $T_0$ for all $x \in X$.  
\end{proposition}

\begin{proof}
Consider the exact sequence
$$
0 \ra U(-n\Delta-X \times \{x\}) \ra U(-n\Delta) \ra U(-n\Delta)|X \times \{x\} \ra 0
$$
on $X \times X$.
Since $U(-n\Delta)|X \times \{x\} \simeq W^* \otimes V(-nx)$ is of degree $-1$ and $W^* \otimes V$ is semistable, this gives
$h^0(W^* \otimes V(-nx)) = 0$ and thus
$$
0 \ra H^1(U(-n\Delta-X\times \{x\}) \ra H^1(U(-n\Delta)) \ra H^1(W^* \otimes V(-nx)).
$$
We claim that
\begin{equation} \label{eq3}
C := h^1(U(-n\Delta)) - h^1(U(-n\Delta - X \times \{x\})) \geq 3.
\end{equation}

According to Lemmas \ref{lem2.2} and \ref{lem2.3},
$$ 
\begin{array}{ll}
C & = gh^0(W^* \otimes V) + l(n-l)(n-1) - h^1(U(-\Delta-X \times \{x\})).
\end{array}
$$
Now the exact sequence 
$$
0 \ra U(-\Delta-X \times \{x\}) \ra U(-X \times \{x\}) \ra U(-X \times \{x\})|\Delta \ra 0
$$
implies 
$$
\begin{array}{ll}
h^1(U(-\Delta-X \times \{x\}) &\leq h^0(U(-X \times \{x\})|\Delta) + h^1(U(-X \times \{x\}))\\ 
&= h^0(W^* \otimes V(-x)) + gh^0(W^* \otimes V).
\end{array}
$$
Hence
$$
C \geq l(n-l)(n-1) - h^0(W^* \otimes V(-x)).
$$

According to the above mentioned theorem of Hirschowitz, either 
$H^0(W^* \otimes V) = 0$ or $H^1(W^* \otimes V) = 0$. 
In the first case also $H^0(W^* \otimes V(-x)) = 0$ and thus
$$
C \geq l(n-l)(n-1) \geq 3.
$$
In the second case Riemann-Roch implies 
$$
h^0(W^* \otimes V(-x)) \leq h^0(W^* \otimes V) = (n+1-g)l(n-l) -1
$$
and thus, for $g \geq 3$,
$$
C \geq l(n-l)(g-2) + 1 \geq 3.
$$
We have thus proved (\ref{eq3}) in all cases.
This implies that the codimension of the union of the kernels of (\ref{eqn3}) for $x \in X$ is at least 2. 
Hence there is a vector subspace $T_0$ of dimension 2 meeting this union in 0 only. 
\end{proof}

\section{The case $n=2$}

\noindent
Now suppose $n=2$, which implies $l=1$. So $V$ and $W$ are line bundles with $\deg(W^* \otimes V) = 1$. 
In this case the proof of Proposition \ref{prop2.4} fails. In fact, we have to choose $V$ and $W$ such that
$$
W^* \otimes V \simeq \cO(x_0)
$$
for some fixed $x_0 \in X$. Then Lemmas \ref{lem2.1} and \ref{lem2.2} remain true and so does Lemma \ref{lem2.3} except when $x=x_0$.

\begin{proposition} \label{prop3.1}
For $n=2$, there is a $(g-1)$-dimensional vector subspace $T_1 \subset H^1(U(-2\Delta))$ such that the restriction map
$$
H^1(U(-2\Delta)) \ra H^1(W^* \otimes V(-2x))
$$
is injective on $T_1$ for all $x \in X$.
\end{proposition}

\begin{proof}
Since $h^0(W^* \otimes V) = 1$, Lemma \ref{lem2.2} says that 
$$
h^1(U(-2\Delta)) = 2g.
$$
Lemma  \ref{lem2.3} implies that, if $x \neq x_0$, then
\begin{equation} \label{eqn4}
h^1(U(-2\Delta - X \times \{x\})) = h^1(U(-\Delta - X \times \{x\})) + g-1.
\end{equation}
If $x=x_0$, then the same proof gives 
\begin{equation} \label{eqn5}
h^1(U(-2\Delta - X \times \{x\})) \leq h^1(U(-\Delta - X \times \{x\})) + g.
\end{equation}
Now consider the exact sequence
\begin{equation} \label{eqn6}
0 \ra U(-\Delta-X \times \{x\}) \ra U(-X \times \{x\}) \ra U(-X \times \{x\})|\Delta \ra 0
\end{equation}
on $X \times X$.
Since under the identification of $\Delta$ with $X$, 
$$
U(-X \times \{x\})|\Delta \simeq \cO(x_0-x),
$$
we get, for $x \neq x_0$,
$$
0 \ra H^1(U(-\Delta - X \times \{x\})) \ra H^1(U(- X \times \{x\})) \stackrel{\varphi}{\ra} H^1(\cO(x_0-x)).
$$
The map $\varphi$ is surjective, since its dual is the canonical injection 
$$
H^0(K_X(x-x_0)) \ra \mbox{Hom}(H^0(\cO(x_0)),H^0(K_X(x))) = H^0(K_X(x)).
$$
Hence 
$$
\begin{array}{ll}
h^1(U(-\Delta - X \times \{x\})) & = h^1(U(-X \times \{x\})) - h^1(\cO(x_0-x))\\
& = h^0(\cO(x_0))h^1(\cO(-x)) - h^1(\cO(x_0-x))\\
& = g - (g-1) = 1.
\end{array}
$$
If $x = x_0$, the map $\varphi$ is still surjective and thus an isomorphism. So (\ref{eqn6}) implies 
$$ 
h^1(U(-\Delta - X \times \{x\})) = h^0(\cO(x_0-x)) = 1.
$$
Now (\ref{eqn4}) and (\ref{eqn5}) give
\begin{equation} \label{eqn7}
h^1(U(-2\Delta-X \times \{x\}))  \left\{ \begin{array}{lll}
                                        \leq g+1 & if & x=x_0,\\
                                        = g& if & x \neq x_0. 
                                         \end{array} \right.
\end{equation}
Now
$$
0 \ra U(-2\Delta-X \times \{x\}) \ra U(-2\Delta) \ra U(-2\Delta)|X \times \{x\} \ra 0
$$ 
gives 
$$
0 \ra H^1(U(-2\Delta-X \times \{x\})) \ra H^1(U(-2\Delta)) \ra H^1(W^* \otimes V(-2x)).
$$
So the kernel of the restriction map is $H^1(U(-2\Delta-X \times \{x\}))$ which, together 
with (\ref{eqn7}), implies the assertion as in the proof of Proposition \ref{prop2.4}.
\end{proof}

\section{Proof of the Theorem for $g \geq 3$}

\noindent
We want to consider extensions of the form
$$
0 \ra q_1^*V(-(n-l)\Delta) \ra E \ra q_1^*W(l\Delta) \ra 0   \eqno(e)
$$
on $X \times X$.  
The extension $(e)$ is classified by an element $e \in H^1(U(-n\Delta))$. The restriction of $(e)$ to $X \times \{x\}$ is the extension
$$
0 \ra V(-(n-l)x) \ra E_x \ra W(lx) \ra 0
$$
corresponding to the image of $e$ in $H^1(W^* \otimes V(-nx))$. We can therefore choose a vector subspace $T_0$ 
of $H^1(U(-n\Delta))$ of dimension 2 such that, for all $0 \neq e \in T_0$, the image of $e$ in $H^1(W^* \otimes V(-nx))$ is non-zero.
Note that 
$$
\begin{array}{ll}
\det E_x & = \det (V(-(n-l)x)) \otimes \det (W(lx))\\
& = \det V \otimes \cO(-l(n-l)x) \otimes \det W \otimes  \cO(l(n-l)x)\\
& = \det V \otimes \det W
\end{array}
$$
for all $x$. On the other hand, by \cite[Lemma 2.1]{ram}, provided $V$ and $W$ are stable, the bundle $E_x$ is stable 
for all $0 \neq e \in T_0$ and all $x \in X$.\\

Let $\PP^1 = P(T_0)$ and consider the product variety $X \times X \times \PP^1$. Let $p_i$ and $p_{ij}$ denote the projections of 
$X \times X \times \PP^1$. 
The non-trivial extensions of the 
form $(e)$ with $e \in T_0$ form a family parametrized by $\PP^1$ which has the form 
(see for example \cite[Lemma 2.4]{ram})
\begin{equation} \label{eqn8}
0 \ra p_1^*V \otimes p_{12}^*\cO(-(n-l)\Delta)) \ra \cE \ra p_1^*W \otimes p_{12}^*\cO(l\Delta) \otimes p_3^*(\tau^*) \ra 0,
\end{equation} 
where $\tau$ is the tautological hyperplane bundle on $\PP^1$. \\

{\it Proof of the Theorem}. By what we have said above, $\cE$ is a family of stable bundles on  $X$ of fixed determinant 
$L = \det V \otimes \det W$ parametrized by $X \times \PP^1$. This gives a morphism
$$
f: X \times \PP^1 \ra M
$$
such that 
$$
(\mbox{id} \times f)^* \cU \simeq \cE \otimes p_{23}^*(N)
$$
for some line bundle $N \in \mbox{Pic}(X \times \PP^1)$. Considering
$$
\cE_x = \cE|\{x\} \times X \times \PP^1
$$
as a bundle on $X \times \PP^1$, we have 
$$
f^*\cU_x \simeq \cE_x \otimes N.
$$
Hence, in order to complete the proof of the theorem, it suffices to show that the bundle $\cE_x \otimes N$ determines the point $x$.

For this we compute the Chern class $c_2(\cE_x \otimes N)$ in the Chow group $\mbox{CH}^2(X \times \PP^1)$.

From (\ref{eqn8}) we get
\begin{equation} \label{eqn9}
c_1(\cE) = p_1^*\beta - (n-l)p_3^*h
\end{equation}
where $\beta$ is the class of $\det V \otimes \det W$ in $\mbox{CH}^1(X)$ and $h$ is the positive generator of $\mbox{CH}^1(\PP^1)$.

For the computation of $c_2(\cE)$ we use the formula 
$$
c_2(\cF \otimes \cL) = c_2(\cF) + (r-1)c_1(\cF)c_1(\cL) + {r \choose 2} c_1(\cL)^2
$$                                                                                                                                       
for any vector bundle $\cF$ of rank $r$ and any line bundle $\cL$.                                                                    

The only terms in $c_2(\cE)$ which can possibly survive in $c_2(\cE_x)$ when restricting are those involving $[\Delta]h$.
So $c_2(p_1^*V \otimes p_{12}^* \cO(-(n-l)\Delta))$ does not contribute. The coefficient of $[\Delta]h$ in 
$c_2(p_1^*W \otimes p_{12}^*\cO(l\Delta) \otimes p_3^*(\tau^*))$ is ${n-l \choose 2}(-2l) $
and the coefficient of $[\Delta]h$ in 
$$
c_1(p_1^*V \otimes p_{12}^*\cO(-(n-l)\Delta))\cdot c_1(p_1^*W \otimes p_{12}^*\cO(l\Delta) \otimes p_3^*(\tau^*))
$$
is $-l(n-l)(-(n-l)) = l(n-l)^2$. This implies 
$$
c_2(\cE_x) = l(n-l)(-(n-l-1) + n-l)(x \times p) = l(n-l)(x \times p),
$$
where $p$ is the class of a point in $\PP^1$. 

Hence, using (\ref{eqn9}), we get that
$$
c_2(\cE_x \otimes N) = l(n-l)(x \times p) + \gamma
$$
with $\gamma \in \mbox{CH}^2(X \times \PP^1)$ independent of $x$. 

If $\cU_x \simeq \cU_y$, then $l(n-l)((x-y) \times p) = 0$ in $\mbox{CH}^2(X \times \PP^1)$. This is equivalent to
$$
l(n-l)(x-y) = 0 \quad \mbox{in} \quad \mbox{CH}^1(X) = \mbox{Pic}(X).
$$
Hence $x-y$ is a point of finite order dividing $l(n-l)$ in $\mbox{Pic}^0(X)$. But there are only finitely many such points in  
$\mbox{Pic}^0(X)$ and any such point has at most 2 representations of the form $x-y$ ( 2 occurs only if $X$ is hyperelliptic). 
So, for general $x \in X$, there is no $y \in X$ such that $x-y$ is of finite order dividing $l(n-l)$ in $\mbox{Pic}^0(X)$.

Now, as stated in the introduction, the natural morphism $X \ra \cM_0,  \; x \mapsto \cU_x$ is \'etale and surjective.
We have now proved that this \'etale morphism has degree 1.
Hence it is an isomorphism, which completes the proof of the theorem. \hfill$\square$

\end{document}